\documentclass[epic,11pt]{article}
\usepackage{amsfonts,amsmath,amscd,amssymb,amsthm}
\usepackage{fullpage}
\usepackage[alphabetic,initials]{amsrefs}
\usepackage[all]{xy}
\usepackage{hyperref} 
\usepackage{mathrsfs}
\input xy%
\usepackage{fullpage}
\def\convf{\hbox{\space \raise-2mm\hbox{$\textstyle      \bigotimes \atop \scriptstyle \omega$} \space}}
\def\0{{\bar 0}}
\def\1{{\bar 1}}

\def\Z{{\mathbb Z}}

\def\B{{\mathcal B}}

\def\Max{{\operatorname{Max}\;}}

\def\Spec  {{\operatorname{Spec  }\;}}

\def\Rad{\operatorname{Rad\;}}

\def\Ker {{\operatorname{Ker}\;}}
\def\Lie {{\operatorname{Lie}\;}}

\def\atyp{{\operatorname{atyp}}}

\newcommand{\ttk}{\mathtt{k}}
\newcommand{\tte}{\mathtt{e}}

\newcommand{\ttT}{\mathtt{T}}

\newcommand{\itema}{\item[{{\rm$($a$)$}}]}
\newcommand{\itemb}{\item[{{\rm$($b$)$}}]}
\newcommand{\itemc}{\item[{{\rm$($c$)$}}]}
\newcommand{\itemd}{\item[{{\rm$($d$)$}}]}
\newcommand{\iteme}{\item[{{\rm$($e$)$}}]}
\newcommand{\itemf}{\item[{{\rm$($f$)$}}]}

\newcommand{\noi}{\noindent}

\newcommand{\ga}{\alpha}
\newcommand{\gb}{\beta}
\newcommand{\gc}{\gamma}

\newcommand{\Gd}{\Delta}

\newcommand{\gd}{\delta}

\newcommand{\gO}{\Omega}

\newcommand{\gt}{\tau}

\newcommand{\gl}{\lambda}

\newcommand{\gep}{\epsilon}

\newcommand{\I}{\mathcal I}

\newcommand{\ot}{\otimes}

\newcommand{\fg}{\mathfrak{g}}\newcommand{\fgl}{\mathfrak{gl}}
\newcommand{\fsl}{\mathfrak{sl}}\newcommand{\osp}{\mathfrak{osp}}

\newcommand{\fh}{\mathfrak{h}}

\newfont{\eufm}{eufm10 scaled\magstep1}

 \newcommand{\ti}{\times}

\newcommand{\bca}{\bigcap}
\newcommand{\bcu}{\bigcup}
\newcommand{\bop}{\bigoplus}

\newcommand{\cO}{\mathcal{O}}

\newcommand{\cI}{\mathcal{I}}

\newcommand{\cS}{\mathcal{S}}
\newcommand{\cT}{\mathbb{T}}

\newcommand{\cV}{\mathcal{V}}

\newcommand{\cW}{\mathfrak{W}}

\newcommand{\ey}{\end{eqnarray}}
\newcommand{\by}{\begin{eqnarray}}
\newcommand{\nn}{\nonumber}

\newcommand{\bco}{\begin{conjecture}}
\newcommand{\ba}{\begin{alg}}
\newcommand{\ea}{\end{alg}}
\newcommand{\eco}{\end{conjecture}}
\newcommand{\bpf}{\begin{proof}}
\newcommand{\epf}{\end{proof}}
\newcommand{\bt}{\begin{theorem}}

\newcommand{\et}{\end{theorem}}

\newcommand{\br}{\begin{rem}}
\newcommand{\er}{\end{rem}}
\newcommand{\brs}{\begin{rems}}
\newcommand{\ers}{\end{rems}}
\newcommand{\bi}{\begin{itemize}}
\newcommand{\ei}{\end{itemize}}
\newcommand{\bl}{\begin{lemma}}
\newcommand{\bsul}{\begin{sublemma}}
\newcommand{\esul}{\end{sublemma}}
\newcommand{\bp}{\begin{proposition}}
\newcommand{\be}{\begin{equation}}
\newcommand{\bc}{\begin{corollary}}
\newcommand{\bexs}{\begin{examples}}
\newcommand{\eexs}{\end{examples}}
\newcommand{\bexa}{\begin{example}}
\newcommand{\eexa}{\end{example}}
\newcommand{\bex}{\begin{exercise}}
\newcommand{\eex}{\end{exercise}}
\newcommand{\btab}{\begin{tab}}
\newcommand{\etab}{\end{tab}}
\newcommand{\el}{\end{lemma}}
\newcommand{\ep}{\end{proposition}}
\newcommand{\ee}{\end{equation}}
\newcommand{\ec}{\end{corollary}}
\newcommand{\Bc}{\begin{center}}
\newcommand{\Ec}{\end{center}}
\newcommand{\bh}{\begin{hyp}}
\newcommand{\eh}{\end{hyp}}
\newcommand{\bhs}{\begin{hyps}}
\newcommand{\ehs}{\end{hyps}}
\newcommand{\bd}{\begin{dfn}}
\newcommand{\ed}{\end{dfn}}

\newcommand{\bn}{\begin{notn}}
\newcommand{\en}{\end{notn}}

\begin{document}
\title{Table of Contents}

\newtheorem*{bend}{Dangerous Bend}

\newtheorem{thm}{Theorem}[section]
\newtheorem{hyp}[thm]{Hypothesis}
 \newtheorem{hyps}[thm]{Hypotheses}
\newtheorem{notn}[thm]{Notation}

  \newtheorem{rems}[thm]{Remarks}

\newtheorem{conjecture}[thm]{Conjecture}
\newtheorem{theorem}[thm]{Theorem}
\newtheorem{theorem a}[thm]{Theorem A}
\newtheorem{example}[thm]{Example}
\newtheorem{examples}[thm]{Examples}
\newtheorem{corollary}[thm]{Corollary}
\newtheorem{rem}[thm]{Remark}
\newtheorem{lemma}[thm]{Lemma}
\newtheorem{sublemma}[thm]{Sublemma}
\newtheorem{cor}[thm]{Corollary}
\newtheorem{proposition}[thm]{Proposition}
\newtheorem{exs}[thm]{Examples}
\newtheorem{ex}[thm]{Example}
\newtheorem{exercise}[thm]{Exercise}
\numberwithin{equation}{section}%
\setcounter{part}{0}
\newcommand{\drar}{\rightarrow}
\newcommand{\lra}{\longrightarrow}
\newcommand{\rra}{\longleftarrow}
\newcommand{\dra}{\Rightarrow}
\newcommand{\dla}{\Leftarrow}
\newcommand{\rl}{\longleftrightarrow}

\newtheorem{Thm}{Main Theorem}


\newtheorem*{thm*}{Theorem}
\newtheorem{lem}[thm]{Lemma}
\newtheorem*{lem*}{Lemma}
\newtheorem*{prop*}{Proposition}
\newtheorem*{cor*}{Corollary}
\newtheorem{dfn}[thm]{Definition}
\newtheorem*{defn*}{Definition}
\newtheorem{notadefn}[thm]{Notation and Definition}
\newtheorem*{notadefn*}{Notation and Definition}
\newtheorem{nota}[thm]{Notation}
\newtheorem*{nota*}{Notation}
\newtheorem{note}[thm]{Remark}
\newtheorem*{note*}{Remark}
\newtheorem*{notes*}{Remarks}
\newtheorem{hypo}[thm]{Hypothesis}
\newtheorem*{ex*}{Example}
\newtheorem{prob}[thm]{Problems}
\newtheorem{conj}[thm]{Conjecture}

\title{Weyl groupoids and superalgebraic sets.}
\author{Ian M. Musson
 \\Department of Mathematical Sciences\\
University of Wisconsin-Milwaukee\\ email: {\tt
musson@uwm.edu}}
\maketitle
\begin{abstract} This paper is a contribution to the study of the geometry of algebras related 
the Weyl groupoid initiated in \cite{M22}. The Nullstellensatz gives a bijection between radical ideals of such an algebra and their zero loci, the superalgebraic sets.  Such sets are exactly the (Zariski) closed sets that are invariant under the action of a suitable groupoid, and the smallest superalgebraic set containing a given closed set can be described explicitly.  Here we give several examples of superalgebraic sets. We also give several 
characterizations of Laurent supersymmetric polynomials. These adapt to 
unite several definitons of one of  the algebras of interest, $J(G)$ that  
may be found in the literature. 
\end{abstract}
\noi 
\section{Background. } \label{Br}  
Throughout  $[k]$ denotes  the set of the first $k$ positive integers and iff means if and only if. We work over an algebraically closed field  $\ttk$  of characteristic zero. In this Section we summarize the results from \cite{M22} that we will need later.  
Let $\ttk$ be an algebraically closed field of characteristic zero.  
\subsection{Algebras related to the Weyl groupoid. } \label{Cls} Let  $\fg= 
\fg(A,\gt)$ be a contragredient Lie superalgebra as defined in  \cite{Kac1} or \cite{M} Equation  (5.2.1).  We adopt the terminology  Kac-Moody Lie superalgebra of finite type or KM algebra 
for this  algebra, which is more descriptive, and becoming more current. 
 Let  $Z(\fg)$ be the center of the enveloping algebra $U(\fg)$.
 There is an injective Harish-Chandra map $\psi:Z(\fg)\lra S(\fh)^W$, whose image  $I(\fh)$ 
we describe in \eqref{ltal} below. Let $\fh$ and $W$  be a Cartan subalgebra and  the Weyl group of the reductive Lie algebra $\fg_0$.
Let $(\;,\;)$  be a non-degenerate, symmetric,   $W$-invariant bilinear form on
$\fh^*$ given by \cite{M}  Theorem 5.4.1 and Remark 5.4.2.  Non-degeneracy  of the form $(\;,\;)$ 
 yields an isomorphism $\fh \lra \fh^*$, denoted  $h_\ga\mapsto \ga$. 
For  $\alpha \in 
\Gd_{iso}$, the set of isotropic roots of $\fg$,  set \be \label{rtal}  \Pi_{\ga}=
 \{ \lambda \in   \fh^* | (\lambda,\ga) =0\}.
\ee 
By work of Gorelik and Kac,  \cite{Gk} \cite{Kac4}
see   also  \cite{M}, Theorem 13.1.1, the image of $\psi$ is 
\be \label{ltal}  I(\fh)= \{f\in S(\mathfrak{h})^W| f(\lambda) = f(\lambda + t \alpha) \mbox{ for all }\alpha \in 
\Gd_{iso}, \gl \in \Pi_{\ga} \mbox{ and } t \in \ttk\}.\ee Several other definitions of $I(\fh)$ can be given based on \cite{M}, Lemma 12.1.1.
\\ \\
The other algebra of interest is the reduced 
Grothendieck $\Z$-algebra 
or supercharacter ring $J(G)$, constructed from finite dimensional representations of $G$. 
Here $G$ is a suitable Lie supergroup with Lie superalgebra $\fg$. 
There is a description of $J(G)$ parallel to \eqref{ltal}. Let 
$\mathbb{T}$ be a maximal torus of $G_0$  and $X(\mathbb{T})$ the character group of $\mathbb{T}$. 
Let  
 $D_\ga$ be the derivation of  $\Z[X(\cT)] $ given by  
\be \label{pvp} D_\ga(\tte^\gb) = (\ga, \gb) \tte^\gb.\ee
Then in the KM case, we have by \cite{SV2} Equation (1), 
\be \label{yta}  J(G)= \{f \in \Z[X(\cT)]^W| D_\ga f \in (\tte^\ga -1)  \mbox{ for all }\alpha \in 
\Gd_{iso}\}.\ee 
Various conditions equivalent to that defining $J(G)$ are given in Section \ref{csd}. To bring out the analogy with $I(\fh)$ we  extend scalars to $\ttk$ and study $J(G) \ot_{\Z}\ttk$.  
The action of 
the continuous Weyl groupoids $\cW^c$  and $\cW_*^c$ on 
$\fh^*$ or $\cT$ is 
defined  in \cite{M22} and it is shown that the algebras $I(\fh)$ and  $J(G) \ot_{\Z}\ttk$  arise as invariant rings for these actions.
All orbits of these groupoids are  closed under these actions.  In this paper we will only be concerned with the case where $G=GL(m|n).$ 
\bexa{\rm If $\fg=\fgl(m|n)$ and  $G=GL(m|n)$ the algebras $I(\fh)$ and  $J(G) \ot_{\Z}\ttk$ are the algebras of 
 supersymmetric polynomials and Laurent supersymmetric polynomials over $\ttk$ respectively.}
\eexa 
\subsection{The  Nullstellensatz.}\label{sn} For simplicity we state some results from \cite{M22} Section 4 in the case of $I(\fh)$ only, but similar results also hold for $ J(G) \ot_{\Z}\ttk$. 
\noi First we recall  the weak Nullstellensatz.
\bt\label{hen}
  If $m$ is a maximal ideal of  $I(\fh)$, then there is a maximal ideal $M$ of  $S(\fh)$ such that $m = M \cap I(\fh)$.\et
\noi For $\gl \in \fh^*,$ let $M_\gl\in \Max S(\fh)$ be the  ideal of functions  vanishing at $\gl$ and $m_\gl =M_\gl\cap I(\fh)$.
\\ \\
The strong Nullstellensatz is deduced from Theorem \ref{hen}. 
\noi If $I$ is a subset of $ I(\fh)$, let
$\cV(I) = \{x\in \fh^*| f(x) =0 \mbox{ for all } f \in I\}$.   
Such a set is called an {\it superalgebraic set}. If instead
$I$ is a subset of $S(\fh)$, we say that $\cV(I)$ is an {\it algebraic set}
(or {\it closed}) in $\fh^*$. 
Thus any superalgebraic set is algebraic.
In addition if
$V$ is a subset of $\fh^*$, 
 set $$\cI_ {I(\fh)}(V) = \{f\in  {I(\fh)}| f(x) =0 \mbox{ for all } x \in V\}.$$ In Section \ref{csc},  we will also need
$$\cI(V) = \{f\in S(\fh)| f(x) =0 \mbox{ for all } x \in V\}.$$
\bt \label{boa}  The maps $\cI_ {I(\fh)}$ and $\cV$ are inverse bijections between the set of superalgebraic sets in $\fh^*$, and the set of radical ideals in $ I(\fh)$. \et
\noi We mention a key ingreditent in the proof of Theorem \ref{hen}.  If $R$ is a  commutative ring  and $T\in  R$, is a non-zero divisor we denote the localization by $R_T$. 
Let $\gO$ be a  $W$-invariant set
of representatives of  the orbits action of $ \{\pm 1\}$ on  $\Gd_{iso}$ and  define 
\be \label{f101} T=\prod_{\beta\in\gO} h_{\beta}.\ee
 We can  choose $ \gO$ so that $T$ is $W$-invariant. 
Then    $TS(\fh)^W \subseteq {I(\fh)}$ and we have an equality of localizations ${I(\fh)}_T = S(\fh)^W_T.$
Then if  $R={I(\fh)}$ or $S(\fh)^W$,  there are order preserving bijections  
\be \label{yok} \{P \in \Spec   R|  T \notin P \} \leftrightarrow \Spec   R_T.\ee
\noi  
This gives a bijection
\be \label{yak} \{P \in \Spec   S(\fh)^W| T \notin P\} \leftrightarrow  \{P \in \Spec   {I(\fh)}| T \notin P\}
,\ee given by $P \mapsto P\cap {I(\fh)}$.  It follows easily that there is  a disjoint union
\be \label{gnu} \Spec   {I(\fh)} = \Spec S(\fh)^W_T \cup \phi^{-1} (\Spec   {I(\fh)}/TS(\fh)^W),  \ee where
$\phi:{I(\fh)} \lra   {I(\fh)}/TS(\fh)^W$ is the natural map.
\bc  \label{pig}  If         $m$ is a maximal ideal of  ${I(\fh)}$ and $T\notin m$, the  maximal ideal $M$ of  $S(\fh)^W$
given by $M= m_T\cap S(\fh)^W$ satisfies  $m = M \cap {I(\fh)}$.   
\ec \noi 
In Corollary \ref{pig}, $m$ and $M$ have the same set of zeroes in $\fh^*$, so by the usual Nullstellensatz  $M  = \Rad mS(\fh)$.  However $M$ can strictly contain $mS(\fh)$, see Example 
\ref{e8}.\\ \\
To study prime, in particular maximal ideals in ${I(\fh)}/TS(\fh)^W$, we use induction and the Duflo-Serganova functor to relate this ring to ${I(\fh_x)}$, where $ \fh_x$ is a Cartan subalgebra of a Lie superalgebra $\fg_x$ of smaller rank.
\subsection{Orbits of $\cW^c$  on $\fh^*$.} \label{orb} 
An {\it iso-set} is a linearly independent set of mutually orhogonal isotropic roots. Define the  {\it degree of atypicality} of $\gl\in \fh^*$ to be  
$$\atyp \; \gl =  \max\{s|(\gl,A)=0 \mbox{ for some iso-set } A \mbox{ with } |A|  = s \} 
.$$ 
\noi 
We say $\gl $ is typical if  $\atyp \; \gl =  0.$
\bt\label{sgs}   Let $F(\gl) $  be an iso-set   with $|F(\gl)| =\atyp\;\gl$. Then 
\bi 
\itema  $\cW^c{\gl} =\bcu_{w\in W} w(\gl +\sum_{\ga\in F(\gl)} \ttk\ga) $, 
\itemb $\dim \cW^c  \gl=  \atyp\;\gl$ and every $\cW^c$-orbit 
is closed.  
\ei\et
\subsection{The Geometry of Superalgebraic Sets.}\label{csc} For simplicity we assume
$\fg=   \fsl(m|n)$ with $m\neq n$ or $\fgl(m|n)$. 
Iidentify $\fh^*$  with 
$\ttk^{m|n}$  and use coordinate functions $x_1,\ldots , x_{m}, y_1,\ldots ,y_n$. Write an element of $\ttk^{m|n}$   as an $m +n$-tuple
\be \label{ear} \gl = (a_1,\ldots, a_{m}, b_1,\ldots ,b_n) =  \sum_{k=1}^{m}  a_k \gep_k + \sum_{\ell=1}^{n}  b_\ell \gd_\ell, 
\ee
where $a_i = x_i(\gl)$ and $b_i = y_i(\gl)$. 
We have $\sum_{i=1}^{m}  a_i = \sum_{i=1}^{n}  b_i$  if  $\fg=\fsl(m|n)$.
We define the S-topology on $\fh^*$ by declaring that the S-closed sets are the $\cW^c$-stable  (Zariski) closed sets. We review a procedure which gives  the $S$-closure $V^S$ of an arbitrary closed set $V=V_0=\cV(I)$ with $I$ a radical ideal in $S(\fh)$.  
 If the defining equations for $V$ are known, this can be done algorithmically using Groebner bases, see Example \ref{qhs}. We can  assume  that $V$, and hence also  $I$ are  $W$-invariant. 
Let 
$$\atyp \;V = \max\{s| \mbox{ for some } \gl\in V, \atyp \; \gl = s \}.$$
Suppose $r=\atyp \;V$ and choose an iso-set  $S=\{\ga_1, \ga_2, \ldots, \ga_r\}$ with  $(\gl,S)=0$ for some $\gl\in V$. Without loss of generality we can assume  $\ga_i = \gd_{n+1-i} - \gep_{m+1-i}$ for all $i\in [r]$.
Set  $A(i)=\{\ga_1,\ldots,\ga_i \} $ and 
$V_{A(i)}= V\cap \bca_{\ga\in {A(i)}} \Pi_\ga$. 
Write $\gl\in \ttk^{m|n}$, in the form \eqref{ear}.
Then define   
\be \label{e6} p_{{A(i)}}(\gl) = 
\sum_{\ell=1}^{m-i}   a_{\ell}  \gep_\ell + \sum_{\ell=1}^{n-i}b_\ell \gd_\ell \ee  
\noi  and 
\by \label{e5}
K_{{A(i)}}&=&  \cI(V_{A(i)}) \cap \ttk[x_1,\ldots ,x_{m-i}, y_1,\ldots,  \ldots , y_{n-i}]. \ey
The ideal $K_{{A(i)}}$ is called an {\it elimination ideal}, \cite{CLO} Chapter 3.
 We denote the closure of $U\subseteq \fh^*$ by 
$\overline{U}$. If $ U_{{A(i)}}= p_{\ga}(V_{{A(i)}})$, then by \cite{CLO} Theorem 3.2.3, $\overline{U}_{{A(i)}}=\cV(K_{A(i)})$.
Let $ L_{A(i)}$ be the ideal of $S(\fh)$ generated by $K_{A(i)}$ and $h_{\ga_1},\ldots, h_{\ga_i}.$ 
For all $ i\in  [r]$, $V^S$ contains the sets 
$$ Y_{{A(i)}}:=\overline{U}_{{A(i)}} +\sum_{j=1}^i\ttk  {\ga_j},$$ and we have 
$Y_{{A(i)}}  =\cV(L_{{A(i)}})$.
However the sets $Y_{{A(i)}}$ are not in general $W$-invariant, so we define $V_i = \bcu_{w\in W} wY_{A(i)}$ and
$I_i = \bca_{w\in W} wL_{A(i)}$. Then 
$V_i = \cV(I_i)$. 
\bt \label{xya} 
The S-closure of $V$ is 
$V^S =   \bcu_{i=0}^r V_i.$  The set $V^S$ is superalgebraic. 
\et\noi  A slightly more involved analysis shows that the superalgebraic sets are exactly the closed sets that are invariant under the Weyl groupoid  $\cW^c$.

\section{Examples.} \label{exs}In all cases $\fg=   \fsl(m|n)$ with $m\neq n$ or $\fgl(m|n)$. 
In the latter case, by \cite{M} Theorems 13.4.1 and  12.4.1, $I(\fh)$ is isomorphic to the algebra of supersymmetric polynomials 
over $\ttk$, which is generated by the  polynomials
 \be \label{prd} p^{(r)}_{m,n} = (x^r_1 + \ldots + x^r_m) + (-1)^{r-1}(y^r_1 +
 \ldots + y^r_n) \ee with $r \geq 1$. If $M \in \Max \B$, and $m = M\cap  I(\fh)$, then 
$$\ttk \subseteq  I(\fh)/m \cong ( I(\fh)+M)/M \subseteq \B/M =\ttk.$$
Thus $m\in \Max  I(\fh)$ and $ I(\fh)/m =\ttk$. Now by Theorem  \ref{hen}, $m = m_\gl$ for some $\gl\in \fh^*$.  So if $q_r =p^{(r)}_{m,n} - p^{(r)}_{m,n}(\gl)$, then the polynomials $q_r$, $r\ge 1$ generate $m$.  By Theorem \ref{sgs} 
$\cV(m_\gl)=\cW^c{\gl}$, so $q_r$ is independent of the choice of $\gl$ in the 
$\cW^c$-orbit. It would be interesting to know when $m_\gl$ is finitely generated. The algebra $ I(\fh)=\bop_{n \ge 0} I(\fh)_n$ is graded, where   $ I(\fh)_n$ consists of all homogeneous polynomials in $I(\fh)$ of degree $n$. The ideal of $I(\fh)$
generated by all   polynomials
$p^{(r)}_{m,n}$, with $r\ge 1$ is the irrelevant ideal $m=\bop_{n \ge 1} I(\fh)_n$.  Since $I(\fh)$ is not Noetherian by  \cite{M} Corollary 13.2.11, $m$ is not finitely generated by \cite{E} Exercise 1.4. Because $m$ is maximal and all  
$p^{(r)}_{m,n}$ vanish at $0\in \fh^*$, Theorem \ref{sgs} implies that $\cV(m)$ is the 
orbit $\cO_0$ of 0 under $\cW^c $. However it is possible for   $\cO_0$  to be cut out by finitely many polynomials in $I(\fh)$ as in the first Example. Thus the radical of a finitely generated ideal need not be finitely generated.  Note $T=\prod_{i\in [m], j\in [n]}(x_i + y_j)\in I(\fh)$ since $T$ 
is supersymmetric. 
We take this as the element $T$ from \eqref{f101}.

\bexa \label{e101} {\rm Let $\fg =\fgl( 2,1)$. 
Set $p_r = x^r +y^r + (-1)^{r-1} z^r$, where $x, y ,z$ are coordinates on $\ttk^{2|1}.$  Then $2T = p_1^2 -p_2$.  Since all points in $\cV(T)$ are atypical, it follows that $\cO_0 =\cV(T,p_1)$.  This orbit is the union of two  lines $(a,0|-a) \cup (0,a|-a) $ as  $a$ runs over $\ttk$.}\eexa

\bexa \label{e8} {\rm Let $\fg =\fgl( 2,1)$. In the situation of Corollary \ref{pig}, we give an example where
$mS(\fh)$ is strictly contained in $M$, that is $mS(\fh)$ is not a radical ideal.
Let $M$ be the ideal of $S(\fh)$ corresponding to the typical point $\gl=(1,1|0)$. Then $m = M\cap  I(\fh)$ is generated by the polynomials
$p_r - p_r(\gl) = p_r-2$  for $r\geq 1$, and $mS(\fh)$ is the ideal of $S(\fh)$ generated by these polynomials. To establish the result we show that $S(\fh)/mS(\fh) \cong \ttk[y]/(y-1)^2$. From now on we work in the ring $S(\fh)/mS(\fh)$, using the same symbol for an element of $S(\fh)$ as its residue in this ring. First we solve
$p_1 =2$ to obtain
\be \label{e1} z = 2-x-y.\ee
\bl \label{e0}In $S(\fh)/mS(\fh)$ we have $x+y =2, z =0, xy = 1$ and $(u-1)^2 = 0$ for $u = x, y$.\el
\bpf This follows from the analysis of the equations $p_r = 2$ for $r =2$ and 3. First substitute \eqref{e1} into $p_r=2$ for $r =2, 3$ and simplify to obtain
\be \label{e2} xy -2x-2y +3 = 0\ee and
\be \label{e3} 2 + 2x^2 + 2y^2 -4x-4y - x^2y -xy^2 +4xy.\ee
Now multiply \eqref{e2} by $x$ and by $y$ and add both resulting equations to
\eqref{e3}. This yields
$x+y =2$, and so $z =0$. Then $ xy = 1$ follows from \eqref{e2} and combining this with $x+y =2$ yields
$(y-1)^2 = 0$. \epf \noi The Lemma shows that $S(\fh)/mS(\fh)$ is a factor algebra of $\ttk[y]/(y-1)^2$.  It remains to show that the remaining equations
$p_r= 2$, equivalently
\be \label{e4} x^r+ y^r = 2\ee
 for $r\geq 3$ impose no further relations on $S(\fh)/mS(\fh)$. From the Binomial Theorem  $u^r \equiv 1 +r(u -1)$ mod $(u-1)^2$ for $u = x, y$.  Thus \eqref{e4} reduces to $x+y = 2.$ The proof shows the following
}\eexa
\bc The ideal $m$ is generated by $p_r-2$ for $r =1, 2, 3.$ \ec
\bexa \label{qxs}{\rm 
 Let $m=n=2$ and consider the $\cW^c$-orbit $\cO$ of $(0,b,0,c)$ where $b,  c$ are fixed and  $b+c\neq 0$.  It is a  union of 4  lines, $L_i$. 
\by
 L_1 &=& \{(a,b|-a,c)\}, \quad  \quad 
 L_2 = \{(b,a|-a,c)\},\nn\\
 L_3 &=& \{(a,b|c,-a)\}, \quad  \quad 
 L_4 = \{(b,a|c-a,)\}. \nn\\
\ey
\noi 
Here $a$ runs over $\ttk$. Note that the $L_i$ form a single $W$ orbit.
}\eexa
\bl \label{cxn} Suppose $b=c=1$, and set 
$q_r=p^{(r)} + (-1)^{r}-1 $  for $r\geq 1$. If $Q$ is the ideal of $I(\fh)$ given by $Q=(T, q_1, q_2).$  Then $\cV(Q) = \cO$.
\el
\bpf Left  to the reader.
\epf
\bexa \label{qds}{\rm For $\fgl(2,2)$, apart from the orbits in Example \ref{qxs}, there is only one other infinite orbit of $\cW^c$ on $\fh^*$. The orbit 
$ \cO_0$ 
 of $(0,0|0,0)$  is the   union of 2 planes,  
\be
 P_1 = \{(a,b|-a,-b)\} \quad \mbox{and} \quad 
 P_2 = \{(b,a|-a,-b)\} .\\
\ee
Here $a, b$ run over $\ttk$.
}\eexa

\bexa \label{nsq1}  {\rm If  $\fg=\fsl(2,1),$ then $\Gd_{iso}=\{\pm\gb, \pm\gc \} $ and $\fh^*=\ttk^{2|1}$. We can choose coordinates so that $h_\gb, h_\gc$ vanish on the horizontal and vertical axes.  We have $I(\fh)=\ttk+h_\gb h_\gc S(\fh)^W$ and $W$ is generated by  a reflection that interchanges the two axes, \cite{M} Example 13.1.3.    
If  $\fg=\osp(3,2),$ then all of the above holds except that  $W$ is generated by  two commuting reflections, each of which interchanges the two axes.  Let $V$ be a circle with center at the origin. Then the $S$-closure of $V$ is the union of $V$ and the two axes.}\eexa

\bexa \label{nsq} The S-closure of a  non-singular quadric. {\rm Let  $\fg=\fgl(n,n)$ and consider the $W$-stable  quadric hypersurface $V=V_0 = \cV(f)\subset \ttk^{m|n}$ defined by $$f_{n}=p^{(2)}_{n,n}-1 = \sum_{i=1}^n x^2_i -  \sum_{i=1}^n y^2_i-1.$$ The S-closure $V^S$ of $V$ involves quadrics in lower dimensional spaces. We compute $V^S$ using the method and notation of Section \ref{csc}. Note  $V$  contains  the typical point  $(1,1,\ldots, 1|0,1,\ldots, 1)$, and $ \atyp \;V =n-1$. Let $\ga_i= \gd_{n+1-i} -\gep_{n+1-i}$.   Then set 
 $A(i)=\{\ga_1,\ldots,\ga_i \} $.  If $\ga = \ga_1$, then $K_\ga$ is the ideal generated by $f_{n-1}$, 
$Y_{\ga}$ is the closure of the set $\{p_{\ga}(\gl) + \ttk \ga|\gl\in V_{\ga}\}$, and $V _{1} =\cup _{w\in W}wY _{\ga}.$ We have $V^S =V\cup V_{1}^S$. 
  More generally $K_{A(i)}$ is generated by $f_{n-i}$,  $Y_{A(i)}$ is the closure of the set $\{p_{A(i)}(\gl) + \sum_{j=1}^i \ttk \ga_j|\gl\in V_{A(i)}\}$, and if $V _{i} =\bcu _{w\in W}wY _{A(i)},$ we have $V^S =\bcu _{i=0}^{n-1} V_{i}$.
} \eexa

\bexa \label{qes}{\rm Let  $\fg=\fgl(2,2)$ and consider the  supersymmetric polynomial  $p = x_1 +x _2 + y_ 1+ y_2$. Then $V = \cV(p)$ contains  typical points, for example $(2,3|-1,-4)$. In addition $V$ contains the orbit 
$ \cO$ 
 of $(0,0|0,0)$  as in Example \ref{qds}, all of whose points have atypicality degree 2. However, it is easy to see that 
$V$ contains no points of atypicality degree 1. Note that if $V $ is   irreducible as a superalgebraic set, but $ V $ has 2 irreducible components as an algebraic set.
}\eexa

\bexa \label{qhs}{\rm 
Consider the polynomials
$$g_1 = x^2 +y+z-1, \quad g_2 =  y^2-y-x^2+x, \quad g_3 = 2x^2y+x^4-x^2, \quad g_4= x^6-4x^4+4x^3 -x^2.$$ These are the polynomials in Equation (3) of \cite{CLO}, Section 3.1, with $x$ and $z$ interchanged.  Let $\I$ be the ideal of $\ttk[x, y,z] $ generated by these polynomials. By \cite{CLO}, page 114 we have $\I_1: = \I \cap \ttk[x, y] = (g_2, g_3, g_4)$ 
and $\I_2: = \I \cap \ttk[x] = (g_4)$. Now we turn this into an example for $\fg=\fgl(3,3)$. 
Let 
 \be Y_0  = \{(x,y,z \mid u, v, w)| g_i(x,y,z)=0 \mbox{ for } i = 1,2,3,4\}.
\ee Since $\cV(\I)$ is non-empty, and there is no restriction on $u, v, w$, $V$ contains points with any degree of atypicality between 0 and 3. For $i=1,2,3$ set $\ga_i = \gd_{4-i } -\gep _{4-i}$ and $A(i)=\{\ga_1,\ldots,\ga_i \} $.
\\ \\
Using the notation of Subsection \ref{csc}, we have
$$Y_{A(1)} = \{(x,y,z \mid u, v, -z)| g_2 \equiv 
  g_3  \equiv g_4 \equiv 0
 \},$$
$$Y_{A(2)} = \{(x,y,z \mid u, -y, -z)|g_4 \equiv 0
 \},$$
$$Y_{A(3)}= \{(x,y,z \mid -x, -y, -z) 
 \}.$$ 
Now let $V_0 =\cup _{w\in W} wY_0$ and $V_i =\cup _{w\in W} wY_{A(i)}$  { for }$ i = 1,2,3.$ Then $V_0$ is $W$-invariant, and the S-closure of $V_0 $ is  $\bcu_{i=0}^3 V_i$.
}\eexa

\section{ On  Laurent supersymmetric polynomials.}\label{csd}
 Let $m,n$ be positive integers.  The ring 
\be \label{crd}
\Lambda_{m,n}=\{f\in \Bbb Z[x_1^{\pm1},\dots,x_m^{\pm1},y_1^{\pm1},\dots,y_n^{\pm1}]^{W}\mid x_i\frac{\partial f}{\partial x_i}+y_j\frac{\partial f}{\partial y_j}\in (x_i-y_j)\}
\ee
 is  called the {\it  ring of Laurent supersymmetric polynomials}, see \cite{Se} Definition 4.1. Let $G$ be the supergroup $GL(m|n)$, and $\fg = \Lie GL(m|n) = \fgl(m,n)$. In \eqref{crd} $W=\cS_m\ti \cS_n$ is the Weyl group of $\fg_0$, a direct product of symmetric groups. 
By \cite{SV2} Proposition 7.2, 
$\Lambda_{m,n}\cong J(G)$ the supercharacter ring of the category of finite dimensional $G$-modules.
We mention some alternate characterizations of Laurent supersymmetric
 polynomials with coefficients. 
If $\ga$ is the isotropic root $\ga = \epsilon_{1} -\delta_{1}$, then 
Lemma  \ref{cod} (e) is a condition on $\ga$ alone without regard to $W$-invariance, so by making appropriate changes, which we leave to the reader, we obtain 
alternate characterizations of other  supercharacter rings $ J(G)$ from 
\cite{SV2}, Section 7. Note condition (e) appears in \eqref{yta}, while (d) appears in 
\cite{HR} Equation (2-3) and (f) is close to conditions  in 
\cite{Gor2} Section 4.2.
\\ \\
Let $\{\epsilon_{1} , \ldots , \epsilon_{m}, \delta_{1} , \ldots , \delta_{n}\}$ be the standard basis for the group $P$ of integral weights of $\fg_0$.  There is an   isomorphism  $P\lra X(\mathbb{T})$,  which we write as $\ga\lra \tte^{\ga}$.
Then 
\[ S= \Z[X(T)] = \Z[x_1^{\pm1},\dots,x_m^{\pm1},y_1^{\pm1},\dots,y_n^{\pm1}],\] where $x_i = \tte^{\epsilon_{i}}$ and $y_i =\tte^{\delta_{i}}.$ 
The direct product of symmetric groups $W = S_m\ti S_n$ acts on $S$ in the obvious way.
  Obviously if $f \in S$ is $W$-invariant, and

\be \label{ape} x_1\frac{\partial f}{\partial x_1}+y_1\frac{\partial f}{\partial y_1}\in (x_1-y_1),\ee
then $f $ is
 Laurent supersymmetric.  We consider the condition (\ref{ape}) independently of $W$-invariance.  For $f\in S^W$, Lemma \ref{cod} shows that condition (\ref{ape})  is equivalent to $f$ being Laurent supersymmetric. 
\\ \\
Let $\ttT, \ttT_1 $ be the subtori of $\mathbb{T}$ defined by 
$\ttT =\Ker x_1 y_1^{-1}, \;\ttT_1 =\Ker x_1 \cap \Ker y_1,$ and let $\ttT_2 = \{c_t=(t,1, \ldots,1| t,1, \ldots,1)|t\in \ttk^*  \}.$   Then we have a direct product $\ttT =\ttT_1 \ttT_2$.
Next set $$R = \Z[X( \ttT_1)]= \Z[x_2^{\pm1},\dots,x_m^{\pm1},y_2^{\pm1},\dots,y_n^{\pm1}] 
,$$ $x = x_1, y = y_1, z_+=(1 - \frac{x}{y}),$  and $z_-=(1 - \frac{y}{x}).$
  Note that $ (1 - z_+)^{-1}= (1-z_-)$ and hence
$ z_- = z_+ (z_+-1)^{-1}$.
It follows that  $S =  R[x^{\pm1}, y^{\pm1}]= R[x^{\pm1}, z_\pm],$
$Sz_+=Sz_-$, and

\be \label{pup} S=R[x^{\pm1},(1 - z_+)^{\pm1}]= R[x^{\pm1},(1 - z_+), (1 - z_-)].\ee
If $\gl =(x_1, x_2, \ldots, x_m| y_1, y_2, \ldots, y_n)$, $t\neq0\in \ttk$ and $f \in S$, set
$f_t(\gl)=f(c_t\gl).$ 
 it is easy to see that derivation 
$D=D_\ga,$  given in \eqref{pvp} satisfies  $Df = x\partial /\partial x + y\partial /\partial y,$
where the partial derivatives vanish on $R.$ The result below is the Laurent analog of \cite{M} Lemma 12.1.1.
\begin{lemma}\label{cod}
Let $z =z_+$. For $f \in S$ the following conditions are equivalent
\begin{itemize}
\itema $f \in R + Sz$
\itemb $f(x = y = t)$ is independent of $t\neq 0$
\itemc
For $\lambda \in \ttT$, and $c_{t} \in \ttT_2$, we have  $f(\lambda) = f(c_{t}\lambda)$.
\itemd $x\partial f/\partial x +y\partial f/\partial y \in (x - y)$
\iteme 
$D_\ga f \in (\tte^\ga -1) $
\itemf If $f = \sum_{\nu \in P} m_\nu \tte^\nu$, then $(\nu,\ga)
\neq 0$ implies that $\sum_{i\in \Z}m_{\nu + i\ga}=0.$ \end{itemize}
\end{lemma}
\noindent \bpf Since $\ga = \epsilon_{1} -\delta_{1}$, we have $\tte^\ga = \frac{x}{y}$.  Thus we have an equality of ideals of $S;
(x - y) =(\tte^\ga -1) 
$, and the equivalence of (d) and (e) follows from the remarks preceding the Lemma.

$(a) \Rightarrow (b)$ If $f \in R $ then $f(x = y =
t) = f$ is independent of $t$, while
if $f \in Sz$ then$f(x = y = t) = 0.$

$(a) \Rightarrow (d)$  This is similar to the proof of $(a)
\Rightarrow (b).$

$(c) \Leftrightarrow (d)$ 
This holds since  $Df $  vanishes on
$\ttT $ iff $f_t({\lambda})$ is constant
for all $\lambda
\in \ttT .$

$(b) \Rightarrow (a)$ Using \eqref{pup}, we can write $f$ uniquely as a finite sum
\begin{equation} \label{fly}
f = \sum_{i \in \Z,j \in \Z}r_{i, j}x^i(1-{z})^j,
\end{equation}
 with $r_{i, j} \in R$ for 
all $i,j.$ Since $(1-z)(x = y = t) = 1$, $f(x = y = t) = \sum_{i, j \in \Z}r_{i, j}t^i,$ so if
(b) holds then for all $i\neq 0$, we have
\begin{equation} \label{emu} r_{i, 0} =-\sum_{j\neq 0} r_{i, j}.\ee
Thus 
\by \label{e9}
f=\sum_{j \in \Z}r_{0, j}(1-{z})^j  &+& + \sum_{i,j \in \Z, i\neq0 }r_{i, j}x^i(1-{z})^j  \nn\\
=\sum_{j \in \Z}r_{0, j}(1-{z})^j  &+& \sum_{i,j \in \Z, i\neq0 \neq j }r_{i, j}x^i((1-{z})^j -1).
\ey
In \eqref{e9} the second sum is contained in $R + Sz$ and we can write the first sum as
$$r_{0, 0}  +
\sum_{j >0}r_{0, j}(1-{z})^j 
+\sum_{j <0}r_{0, j}(1-{z_{-}})^{-j} \in R + Sz. 
$$ 

$(c) \Rightarrow (a)$
Given $\lambda \in \ttT $,
we can find $\lambda' \in \ttT_1 $ and $c_s \in \ttT_2$ such that $\lambda =
c_s\lambda'.$ Hence (c) is equivalent to $f(\lambda) =
f(c_{t}\lambda)$ for all $c_{t} \in \ttT_2$ and $\lambda \in \ttT_1 .$
Write  $f$ as in (\ref{fly}). Then an argument similar to the one used in the previous step of the proof shows that
 if (c) holds, then again \eqref{emu} holds (with $r_{i, j}(\gl)$ replacing $r_{i, j}$) for all non-zero $i\in \Z.$ Since this holds for all  $\lambda \in \ttT_1 ,$ and   $r_{i, j}  \in R = \Z[X( \ttT_1)],$ we see that \eqref{emu} holds exactly as stated.  Thus (a) follows as before. 

$(e) \Leftrightarrow (f)$ .  Define an equivalance relation $\sim$ on $P$ by $\gl\sim\mu$ iff $\gl-\mu\in \Z\ga.$  Choose a representative $\nu  $ for each  equivalence class $[\nu ]$.  Then 

$$f = \sum_{[\nu] } \sum_{i\in \Z}
m_{\nu + i\ga}
\tte^{\nu + i\ga}, \mbox{ 
and } D_\ga(f) = \sum_{[\nu] } \sum_{i\in \Z}
m_{\nu + i\ga} (\nu,\ga)
\tte^{\nu + i\ga}.$$
On the other hand, suppose that 
$g = \sum_{[\nu] } \sum_{i=p(\nu)}^{q(\nu)} 
a_{\nu + i\ga}
\tte^{\nu + i\ga},$ with $p(\nu) \le {q(\nu)}$.  Then $(\tte^\ga -1)g = \sum_{[\nu] } \sum_{i=p(\nu)
}^{q(\nu)} 
(a_{\nu + (i-1)\ga} -a_{\nu + i\ga})
\tte^{\nu + i\ga}.$  Thus 
$D_\ga(f) = (\tte^\ga -1)g $ holds, then for all 
$ \nu$ with $(\nu,\ga)
\neq 0$, we have 
\be \label{grb}
 (\nu,\ga)
m_{\nu + i\ga} = 
a_{\nu + (i-1)\ga} -a_{\nu + i\ga}
\mbox{ whenever } p(\nu) \le i\le {q(\nu)+1}.\ee
If \eqref{grb} holds, then since $a_{\nu + (p-1)\ga}=a_{\nu + (q+1)\ga}=0,$ 
we obtain $\sum_{i\in \Z}m_{\nu + i\ga}=0.$ Conversely if $\sum_{i\in \Z}m_{\nu + i\ga}=0$ holds
for all 
$ \nu$ with $(\nu,\ga)
\neq 0$, then we can solve Equations \eqref{grb}  recursively for the $a_{\nu + i\ga}$ to find $g$ with 
$D_\ga(f) = (\tte^\ga -1)g.$
  \epf





\begin{bibdiv}
\begin{biblist}

\bib{CLO}{book}{
   author={Cox, David},
   author={Little, John},
   author={O'Shea, Donal},
   title={Ideals, varieties, and algorithms},
   series={Undergraduate Texts in Mathematics},
   edition={2},
   note={An introduction to computational algebraic geometry and commutative
   algebra},
   publisher={Springer-Verlag, New York},
   date={1997},
   pages={xiv+536},
   isbn={0-387-94680-2},
}

\bib{DS}{article}{ author={Duflo, Michel}, author={Serganova, Vera V.}, title={On associated variety for Lie
superalgebras}, journal={arXiv:math/0507198.}, date={2005}}

\bib{E}{book}{
   author={Eisenbud, David},
   title={Commutative algebra},
   series={Graduate Texts in Mathematics},
   volume={150},
   note={With a view toward algebraic geometry},
   publisher={Springer-Verlag, New York},
   date={1995},
   pages={xvi+785},
   isbn={0-387-94268-8},
   isbn={0-387-94269-6},
}
\bib{Gk}{article}{ author={Gorelik, Maria}, title={The Kac construction of the centre of $U(\germ g)$ for Lie superalgebras}, journal={J. Nonlinear Math. Phys.}, volume={11}, date={2004}, number={3}, pages={325--349}, issn={1402-9251}, 
}


\bib{Gq}{article}{
   author={Gorelik, Maria},
   title={Shapovalov determinants of $Q$-type Lie superalgebras},
   journal={IMRP Int. Math. Res. Pap.},
   date={2006},
   pages={Art. ID 96895, 71},
   issn={1687-3017},
}

\bib{Gor1}{article}{ author={Gorelik, Maria}, title={Depths and cores in the light of DS-functors.}, journal={arXiv:2010.05721.}, date={2020}}


\bib{Gor2}{article}{ author={Gorelik, Maria}, title={On the character ring of a quasireductive Lie superalgebra.}, journal={arXiv:2203.04932.}, date={2022}}


\bib{GHSS}{article}{ author={  Gorelik, M.}, author={Hoyt, C. }, author={Serganova, V.}, author={ Sherman, A.}, title={The Duflo-Serganova functor, vingt ans après.}, journal={arXiv:2203.00529.}, date={2022}}


\bib{HR}{article}{
   author={Hoyt, Crystal},
   author={Reif, Shifra},
   title={Grothendieck rings for Lie superalgebras and the Duflo--Serganova
   functor},
   journal={Algebra Number Theory},
   volume={12},
   date={2018},
   number={9},
   pages={2167--2184},
   issn={1937-0652},
}


\bib{Kac1}{article}{ author={Kac, V. G.}, title={Lie
superalgebras}, journal={Advances in Math.}, volume={26},
date={1977}, number={1}, pages={8--96}, issn={0001-8708},
}


\bib{Kac4}{article}{author={Kac, V. G.}, title={Laplace operators of infinite-dimensional Lie algebras and theta functions}, journal={Proc. Nat. Acad. Sci. U.S.A.}, volume={81}, date={1984}, number={2, Phys. Sci.}, pages={645--647}, issn={0027-8424}, 
}

\bib{M}{book}{author={Musson, I.M.}, title={Lie Superalgebras and Enveloping Algebras},
   series={Graduate Studies in Mathematics},
   volume={131},
publisher={American Mathematical Society}, place={Providence, RI}, date ={2012}}

\bib{M22}{article}{ author={Musson, I.M.}, title={On the geometry of algebras related to the Weyl groupoid.}, 
journal={In preparation.}, 
date={2022}}


\bib{Se}{article}{
   author={Sergeev, A. N.},
   title={On rings of supersymmetric polynomials},
   journal={J. Algebra},
   volume={517},
   date={2019},
   pages={336--364},
   issn={0021-8693},
}

\bib{SV2}{article}{
   author={Sergeev, Alexander N.},
   author={Veselov, Alexander P.},
   title={Grothendieck rings of basic classical Lie superalgebras},
   journal={Ann. of Math. (2)},
 volume={173},
   date={2011},
   number={2},
   pages={663--703},
   issn={0003-486X},
}
\end{biblist}
\end{bibdiv}
\end{document}